\documentclass[10pt,fleqn]{article}

\usepackage{latexsym,amssymb,amsmath}
\usepackage{lastpage}
\usepackage{amsthm}

\newtheorem{theorem}{Theorem}[section]
\newtheorem*{theoremMain}{Theorem~\ref{thm:cutdown}}
\newtheorem{lemma}[theorem]{Lemma}
\newtheorem{corollary}[theorem]{Corollary}
\newtheorem{conjecture}[theorem]{Conjecture}
\theoremstyle{definition}

\theoremstyle{remark}
\newtheorem{remark}[theorem]{Remark}

\title{New $2$--critical sets in the abelian 2--group}
\author{Carlo H\"{a}m\"{a}l\"{a}inen \\
Centre for Discrete Mathematics and Computing \\
Department of Mathematics \\
The University of Queensland \\
Queensland 4072, Australia \\
carloh@maths.uq.edu.au }
\begin{document}
\maketitle

\begin{abstract}
In this paper we determine a class of critical sets
in the abelian \mbox{2--group} that may be obtained from a greedy algorithm. 
These new critical sets are all 2--critical (each entry intersects
an intercalate, a trade of size $4$) and completes in a top down manner.
\end{abstract}

\setcounter{page}{1}
\renewcommand{\thepage}{\roman{page}}

\section{Introduction}
\renewcommand{\thepage}{\arabic{page}}
\setcounter{page}{1}

Critical sets are minimal defining sets in latin squares
\cite{MR97i:05017}. Some recent work has investigated the structure and
size of critical sets in the latin square $L_s$ derived from the abelian
2--group of order $2^s$ (\cite{MR97f:05031}, \cite{MR99f:05019}). In this paper we
present a new family of critical sets derived from isotopisms of $L_s$.

Section~\ref{sec:definitions} presents background definitions.
Section~\ref{sec:gcs} has basic properties of greedy critical sets.
Then Section~\ref{sec:gcs-ab2group} develops some properties of greedy
critical sets in $L_s$, and Section~\ref{sec:main} completes the proof
of the main result, which is Theorem~\ref{thm:cutdown}.
The Appendices provide extra examples to aid in the
understanding of the Theorem and also have more detail for the inductive
hypotheses.

\section{Definitions}
\label{sec:definitions}

We begin with some definitions.
Let $N^k_n = \{ nk, nk+1, \ldots, nk+n-1 \}$ for integers $k \geq 0$ and
$n > 0$.
A {\em latin square} $L$ of {\em order $n$} is an $n \times n$ array with
rows indexed by $N^{k}_n$,
columns by $N^{k'}_n$,
and with entries from the set $N^{k''}_n$. Further, each 
$e \in N^{k''}_n$ appears exactly once in each row
and exactly once in
each column. This is equivalent to the usual definition where
$k=k'=k''=0$ but allows more flexibility when discussing subsquares.
A {\em partial latin square} 
is an $n \times n$ array where each entry of $N^{k''}_n$ occurs at most once
in each row and at most once in each column. 

A latin square $L$ may also be
represented as a set of ordered triples, where
$(r,c;e) \in L$ denotes the fact that symbol $e$ appears in the cell at
row $r$, column $c$, of $L$. 
The {\em size} of a partial latin square $P$ is the number of filled cells, denoted
by $\left| P \right| = \left| \{ (r,c;e) \mid (r,c;e) \in P \} \right| $.

A partial latin square $L$ of order $n$ is {\em isotopic} to $L'$ (also of
order $n$) if the rows, columns, and entries of
$L$ can be rearranged to obtain $L'$. Specifically, we say that $L$ is isotopic to $L'$
if there exist permutations $\alpha$, $\beta$, $\gamma$ on 
the row labels, column labels, and symbols (respectively) such that
$L' = \{ ( \alpha r, \beta c; \gamma e ) \mid (r,c;e) \in L \}$. We say that
$(\alpha, \beta, \gamma)$ is an {\em isotopism} from $L$ onto $L'$, 
and we write this as 
$L' = (\alpha, \beta; \gamma) L$. We write
$\alpha$ instead of $(\alpha,\iota,\iota)$ when it is clear from the
context that the
columns and entries are left fixed. 

Given a partial latin square $P$ of order $n$, 
we define the partial latin square $P^r = \{ (i,j;k+nr) \mid (i,j;k) \in P  \}$.
Note that if $P$ has symbols selected from $N^{p}_n$, then 
$P^r$ has symbols from $N^{p+r}_n$.
We use this exponent notation
when recursively constructing larger partial latin squares. For example, suppose that
$A$,
$B$,
$C$, and
$D$ are partial latin squares of order $n$. Then by
\begin{center}
	$P = $~\begin{tabular}{|c|c|}
	\hline $A$ & $B$ \\
	\hline $C$ & $D$ \\
	\hline 
	\end{tabular}
\end{center}
we mean the partial latin square $P$ of order $2n$ where
\begin{equation*}
\begin{split}
	P &= 	\{ (i,j;k) \mid (i,j;k) \in A \}
	\cup 	\{ (i,j+n;k) \mid (i,j;k) \in B \} \\
	&\cup 	\{ (i+n,j;k) \mid (i,j;k) \in C \}
	\cup 	\{ (i+n,j+n;k) \mid (i,j;k) \in D \}
\end{split}
\end{equation*}
Let $P$ and $Q$ be partial latin squares of order $n$. Suppose that
$\alpha$, $\beta$, $\gamma$ are bijections between the row, column, and
symbol sets (respectively) of $P$ and $Q$ such that
\begin{enumerate}
\item $P = (\alpha, \beta, \gamma)Q$.
\item $\alpha$ and $\beta$ are monotone.
\end{enumerate}
Then $P$ and $Q$ are said to be {\em similar}, written $P \approx Q$.
Informally, $P$ and $Q$ are similar when
the rows and columns of $Q$ can be 
relabelled (preserving order) to give $Q'$ such that
$P = (\iota, \iota, \gamma)Q'$.

Given a partial latin square $P$ we can define a binary relation $(P, \ll)$
on the elements of $P$ as follows (see also
\cite{awoca-gcs}). For all $(x,y;z), (r,s;t) \in P$,
$(x,y;z) \ll (r,s;t)$ if and only if
\begin{enumerate}
\item $x<r$, or
\item $x=r$ and $y\leq s$.
\end{enumerate}
We can verify that $(x,y;z) \ll (x,y;z)$ so $\ll$ is reflexive.
If $(x,y;z) \ll (r,s;t)$ and 
$(r,s;t) \ll (x,y;z)$ then $x=r$ and $y=s$, so $\ll$ is antisymmetric. Finally,
suppose that $(x,y;z)\ll(r,s;t)$ and $(r,s;t)\ll(u,v;w)$. If $x < r$ then
$x < u$, so $(x,y;z)\ll(u,v;w)$. On the other hand, if $x=r$ and $r < u$ then
$x < u$, so $(x,y;z)\ll(u,v;w)$ again. Finally, if $x=r$ and $r=u$, then
$y \leq s$ and $s \leq v$ so $y \leq v$, which implies that
$(x,y;z)\ll(u,v;w)$. Hence $\ll$ is transitive, and $(P, \ll)$ is
a weak partial order.

In fact, $(P, \ll)$ is a total order since for any
distinct $(x,y;z),(r,s;t)\in P$, either $x<r$, or $r<x$, or
$r=x$ and $y\leq s$ or $s\leq y$. Given a partial latin square 
$P$ we denote the {\em least element} of $(P, \ll)$ by $(rl_P,cl_P;el_P)$
and the {\em greatest element} by $(rg_P,cg_P;eg_P)$. Since $\ll$ is the
only partial order used in this paper we simply say that $(i,j;k)
\in P$ is the least (greatest) element of $P$.

It is convenient to refer to the set of entries occurring in a particular
row or column of a partial latin square $P$. 
For each row 
$i$ of $P$, define 
$\mathcal{R}^i_P = \{ k \mid \mbox{ there exists $j$ such that } (i,j;k) \in P \}$.
Also, for each 
column $j$ of $P$, we define 
$\mathcal{C}^j_P = \{ k \mid \mbox{ there exists $i$ such that } (i,j;k) \in P \}$.
The {\em shape} of a partial latin square $P$ is the set of filled
cells, defined by $\mathcal{S}_P = \{ (i,j) \mid (i,j;k) \in P \}$.

For some partial latin square $P$ we use the following notation to
specify a subsquare:
\[
\mathcal{Q}^k_{i,j}(P) = \{ (x,y;z) \in P \mid i \leq x < i+k, \; j \leq y < j+k \}
\]
We also use this notation for defining subsquares in a partial latin
square. For example, 
$\mathcal{Q}^k_{i,j}(P) = L$ places the order $k$ latin square $L$ into
$P$ starting with the top--left corner at cell $(i,j)$.

Let $P$ be a partial latin square of order $n$ contained in the latin square $L$. 
Without loss of generality, suppose that the rows and columns are
indexed by $N_n = N^0_n$, and that each entry is from $N_n$.
Let
$R \subseteq N_n$, $C \subseteq N_n$, and $S = R \times C$. For each
$(r,c) \in S$, define
\begin{equation*}
S_{r,c} = 
	\begin{cases}
		\emptyset, & \text{ if $(r,c;e) \in P$ for some $e \in N_n$} \\
		N \setminus (\mathcal{R}^r_P \cup \mathcal{C}^c_P), & \text{ otherwise.}
	\end{cases}
\end{equation*}
Then the {\em
array of alternatives} of $S$ with respect to $P$ and $L$ is 
given by $\mathcal{A}(P, S, L) = \{ (r,c; S_{r,c}) \mid (r,c) \in S \}$. For
clarity we write $\mathcal{A}(P, S, L)_{r,c}$ for $S_{r,c}$. 

We say that $\mathcal{A}(P,S,L)$ is similar to $\mathcal{A}(P', S', L')$ if there are relabellings of the
row names, column names and symbols so that the table for 
$\mathcal{A}(P,S,L)$ is equal to the relabelled table for $\mathcal{A}(P', S', L')$.

A partial latin square $T$ forms a {\em latin trade} in
a latin square of order $n$ if there exists a partial latin square
$T'$, the {\em disjoint mate}, such that:
\begin{enumerate}

\item $T$ and $T'$ are of the same order.

\item $\{ (i,j) \mid (i,j;k) \in T \mbox{ for some symbol $k$} \} 
\newline
= \{ (i,j) \mid (i,j;k') \in T' \mbox{ for some symbol $k'$} \}$

\item For each $(i,j;k) \in T$ and $(i,j;k') \in T'$, $k \neq k'$.

\item For each $i \in N$, $\mathcal{R}^i_{T} = \mathcal{R}^i_{T'}$
and $\mathcal{C}^i_{T} = \mathcal{C}^i_{T'}$.

\end{enumerate}

Informally, Condition 2 says that $T$ and $T'$ have the same shape,
Condition 3 says that they are disjoint, and Condition 4 says that $T$
and $T'$ are {\em row balanced} and {\em column balanced}.

Let $L$ and $L'$ be two disjoint latin square of the same order. Let $T = L \setminus L'$
and $T' = L' \setminus L$. Then $T$ and $T'$ form a latin trade. We assume that
all latin trades are nonempty.
A partial latin square $P$ is {\em uniquely completable} if there is
just one latin square $L$ of the same order as $P$ such that
$P \subseteq L$.

A partial latin square $P$ of order $n$ is {\em strongly completable} if
it is uniquely completable to $L$, there is a sequence of partial
latin squares
$P_0 = P \subset P_1 \subset P_2 \subset \ldots \subset P_m = L$
where $m = n^2 - \left| P \right|$, and for each $P_k$ there exists
$r$, $c$ 
such that
$\left| \mathcal{A}(N_n \times N_n, P_k)_{r,c} \right|=1$. 

A partial
latin square $C \subseteq L$ is a {\em critical set} if
\begin{enumerate}

\item $C$ has unique completion to $L$, and

\item no proper subset of $C$ satisfies 1.
\end{enumerate}

A {\em strong critical set} is a critical set that has strong completion.
We say that a (uniquely completable) partial latin square {\em
extends
top down} if, given that rows $0, 1, \ldots, i$ are filled in,
then row $i+1$ can be shown to have unique extension. If all rows can be
extended in this manner then the critical set has {\em unique completion
top down}.

\begin{lemma}
Let $P$ be a critical set in the latin square $L$ and
$T$ a latin trade in $L$. 
Then $P \cap T \neq \emptyset$.
\end{lemma}

\begin{lemma}
Let $L$ be a latin square and 
$C \subseteq L$ a critical set.
For each $x \in C$ there exists a latin trade $T \subseteq L$ such that 
$C \cap T = \{ x \}$.
\end{lemma}

The latin trade containing the least number of entries is a 
$2 \times 2$ subsquare, known as an {\em intercalate}. Let $C$ be a critical
set in $L$, $c \in C$, and $I \subseteq L$ an intercalate such that
$C \cap I = \{ c \}$. Then $c$ is said to be 
{\em 2--essential}. If all $c \in C$ are
2--essential then $C$ is {\em 2--critical}.

\section{Greedy Critical Sets}
\label{sec:gcs}

Algorithm~A was first presented in \cite{awoca-gcs}. Given a partial
latin square $P$ with unique completion,
and a bijection on its cells, the algorithm produces a critical set.  \\

\begin{figure}
\noindent {\bf Algorithm A} \\
\noindent {\bf Input: } Partial latin square $P$ of order $n$ with
unique completion, and \\
\noindent \phantom{{\bf Input: }} bijection $f : \{1, \ldots, \left| P \right| \} \rightarrow \mathcal{S}_P$.

\begin{tabbing}
\hspace*{0.25in}\=\hspace{3ex}\=\hspace{3ex}\=\hspace{3ex}\=\hspace{3ex}\kill

\> $P_0 \leftarrow P$ \\
\> for $i = 1, \ldots, \left| P \right|$ \\
\> \> let $x$,$y$,$z$ be integers such that $(x,y;z) \in P_{i-1}$ and $f(i) = (x,y)$ \\
\> \> if $P_{i-1} \setminus \{ (x,y;z) \}$ has unique completion then \\
\> \> \> $P_i \leftarrow P_{i-1} \setminus \{ (x,y;z) \}$ \\
\> \> else \\
\> \> \> $P_i \leftarrow P_{i-1}$ \\
\> return $P_{\left| P \right| }$ \\
\end{tabbing}
\end{figure}

\begin{lemma}[Lemma~2.1, \cite{awoca-gcs}]
Let $P$ be a partial latin square that uniquely completes
to $L$. Then for every bijection $f$ over
$\{1, \ldots, \left| P \right| \}$, Algorithm~A returns 
a critical set.
\end{lemma}

\begin{proof}
Algorithm~A works on a sequence of partial latin squares, 
$P_0 = P \supseteq P_1 \supseteq \ldots \supseteq P_m$
where $m = \left| P \right|$.
The initial partial latin square $P_0=P$ has unique completion, and the 
\texttt{if} statement ensures that each $P_i$, for $i > 0$, has unique
completion. Hence $P_m$ has unique completion.

To see that $P_m$ is minimal, suppose otherwise.
Then there is an $x \in P_m$ such that
$P_m \setminus \{x \}$ has unique completion. Also,
let $k$ be the integer such that $f(k) = x$. Then
$P_k$ is the partial latin square 
where $x$ is inspected (and not removed)
by Algorithm~A.
Since $P_m \setminus \{ x \}$ has unique completion, we can add entries
to $P_m \setminus \{ x \}$ until we have precisely $P_k \setminus \{ x \}$. This has
unique completion, yet Algorithm~A apparently did not remove $x$,
a contradiction. Hence $P_m$ is minimal and so $P_m$ is a critical set.
\end{proof}

Since a latin square trivially has unique completion, we get:

\begin{corollary}
\label{cor:alg-A-critical-set}
If the input to Algorithm~A is a latin square of order $n$ then
the output is a critical set for any bijective function $f$.
\end{corollary}

We refer to Algorithm~A as the {\em generalised greedy critical
set} algorithm, and abbreviate this to $\mbox{ggcs}(L, f)$ for
given latin square $L$ and map $f$. 

\begin{lemma}
These two sets are equal:
\[
\{ \mbox{ggcs}(L, f) \mid \mbox{ $f$ is a bijection on 
$\{ 1, \ldots, n^2 \}$
} \}
\]
and
\[
\{ C \mid C \subseteq L \mbox{ and $C$ is a critical set of $L$} \}
\]
for some latin square $L$ of order $n$.
\end{lemma}

Let $f_0 : \{ 1, \ldots, n^2 \} \rightarrow \mathcal{S}_L$
be the bijection defined by
\[
f_0(i) = \left( \left\lfloor \frac{i-1}{n} \right\rfloor , n-i \pmod n \right)
\]
for $1 \leq i \leq n^2$ and $L$ of order $n$. 
Then $f_0$ orders the cells of $L$
from right to left along each row and from the bottom row
to the top row. 
We abbreviate
$\mbox{ggcs}(L, f_0)$
to $\mbox{gcs}(L)$ and call this the 
{\em greedy critical set of $L$}.

We now characterise greedy critical sets in terms of the partial order $\ll$.
Let $L$ be a latin square, and ${\cal I}=\{I \mid I \subset
L\mbox{ and } I\mbox{ is a latin trade}\}$. Each $I\in {\cal I}$ is a
partial latin square implying $I$ has a least element and a greatest
element.

\begin{lemma}[Lemma~2.4, \cite{awoca-gcs}]
Let $C$ be a critical set in $L$. Then $C = \mbox{gcs}(L)$
if and only if for all $(x,y;z)\in C$ there exists an $I\in{\cal I}$
such that $I\cap C=\{(x,y;z)\}$ and $(x,y;z)=(rl_I,cl_I;el_I)$.
\label{lemma:dmd-gcs}
\end{lemma}

\begin{proof}
(if) Algorithm~A with input $L$ and map $f_0$ computes on 
a sequence of partial latin squares
$P_1 = L \supseteq P_2 \supseteq \ldots \supseteq P_m = C$. Suppose that
$f_0(k) = (x,y)$ and that $P_{k-1} \setminus \{ (x,y;z) \}$ completes to 
$L_0=L, L_1, \ldots, L_s$ for some $s \geq 1$. Then there is an $L_i$,
$i \geq 1$, such that $P_{k-1} \cap (L \setminus L_i) = \{ (x,y;z) \}$.
In other words, $T = L \setminus L_i$ is  a latin trade in $L$.
The definition of $f_0$ implies that 
for any $(r,c;e) \in T$ then
either $r > x$, or, if $r=x$ then $c \geq y$. Hence $(x,y;z)$ is the
least element of $T$.

(only if) Assume that for all $(x,y;z)\in C$ there exists an $I\in{\cal I}$
such that $I\cap C=\{(x,y;z)\}$ and $(x,y;z)=(lr_I,lc_I;le_I)$, but $C$
is not the greedy critical set $gcs(L)$.  Let $D= L \cap ((C \setminus
\mbox{gcs}(L)) \cup (\mbox{gcs}(L) \setminus C))$, that is, the intersection
with the symmetric difference.

The set $D$ is a partial latin square and 
has a greatest element
$(gr_D,gc_D;ge_D)$
since $D \neq \emptyset$.
Thus for all $(a,b;c)\in L$ such
that $a>gr_D$, or $a=gr_D$ and $b > gc_D$, 
$(a,b;c)\in C$ if and only
if $(a,b;c)\in gcs(L)$. The reason is that $(a,b;c)$ is not in $D$,
so $(a,b;c) \notin ((C \setminus \mbox{gcs}(L)) \cup (\mbox{gcs}(L)
\setminus C))$.

By the definition of $D$ there are two possibilities:

\begin{enumerate}

\item $(gr_D,gc_D;ge_D)\in C$, and $(gr_D,gc_D;ge_D)\notin gcs(L)$.  
Since
$(gr_D,gc_D;ge_D)$ is in $C$, there exists an $I\in{\cal I}$ such that $I\cap
C=\{(gr_D,gc_D;ge_D)\}$ and $(gr_D,gc_D;ge_D)=(lr_I,lc_I;le_I)$. But for all
$(a,b;c)\in L$ such that $a>gr_D$, or $a=gr_D$ and $b > gc_D$, $(a,b;c)\in
C$ if and only if $(a,b;c)\in gcs(L)$, so $I\cap gcs(L)=\emptyset$, which
is a contradiction.

\item $(gr_D,gc_D;ge_D)\in gcs(L)$, and $(gr_D,gc_D;ge_D)\notin C$. 
Let $k$ be the integer such that 
$f_0(k) = (gr_D,gc_D)$.
Then at step $k$ Algorithm~A removes $(gr_D,gc_D;ge_D)$
and
$P_{k-1} \setminus\{(gr_D,gc_D;ge_D)\}$ is found to have
at least two completions, say $L$ and $L'$. So $T = L \setminus L'$
is a latin trade and the least
element of $T$ is $(gr_D,gc_D;ge_D)$. Once again,
this implies that $T \cap C=\emptyset$, which is
a contradiction.  \end{enumerate}

Hence $D = \emptyset$, which contradicts our original assumption that
$C$ was different to the greedy critical set.
\end{proof}

\begin{corollary}
Let $L$ be a latin square of order $n$ and $G = \mbox{gcs}(L)$. 
If $(i,j;k) \in G$ then $i \neq n-1$ and $j \neq n-1$.
\end{corollary}

\section{Greedy Critical Sets in the Abelian 2--Group}
\label{sec:gcs-ab2group}

We define $L_s$ to be the latin square corresponding to the
abelian 2--group of order $n = 2^s$ and the partial latin square
$P_s \subset L_s$ as in \cite{donovan-2crit}. That is,

\begin{center}
\begin{tabular}{cc}

$P_1 = ~$ \begin{tabular}{|c|c|c|c|c|c|}
\hline 0 & $\;\;$ \\
\hline $\;\;$ & $\;\;$ \\
\hline 
\end{tabular}

&

$L_1 = ~$ \begin{tabular}{|c|c|c|c|c|c|}
\hline 0 & 1 \\
\hline 1 & 0 \\
\hline 
\end{tabular} \\

\end{tabular}
\end{center}
and for $s \geq 2$,
\begin{alignat*}{3}
L_s &= L_1 \times L_{s-1} &= &\{ (x,y;z), (x,y+n/2;z+n/2),
(x+n/2,y;z+n/2),\\
&&	& (x+n/2,y+n/2;z) \mid (x,y;z) \in L_{s-1} \}, \mbox{ and} \\
P_s &= P_1 \otimes P_{s-1} &= &\{ (x,y;z), (u,v+n/2;w+n/2),
(u+n/2,v;w+n/2),\\
	&&	& (u+n/2,v+n/2;w) \mid (u,v;w) \in P_{s-1} \mbox{ and } (x,y;z) \in L_{s-1} \}.
\end{alignat*}
For example, $L_3$ and $P_3$ are:
\begin{center}
\begin{tabular}{cc}
	\begin{tabular}{|c|c|c|c||c|c|c|c|}
	\hline 0&1&2&3&4&5&6&7\\
	\hline 1&0&3&2&5&4&7&6\\
	\hline 2&3&0&1&6&7&4&5\\
	\hline 3&2&1&0&7&6&5&4\\
	\hline
	\hline 4&5&6&7&0&1&2&3\\
	\hline 5&4&7&6&1&0&3&2\\
	\hline 6&7&4&5&2&3&0&1\\
	\hline 7&6&5&4&3&2&1&0\\
	\hline 
	\end{tabular}
&
	\begin{tabular}{|c|c|c|c||c|c|c|c|}
	\hline 0&1&2&3&4&5&6&~\\
	\hline 1&0&3&2&5&4&~&~\\
	\hline 2&3&0&1&6&~&4&~\\
	\hline 3&2&1&0&~&~&~&~\\
	\hline
	\hline 4&5&6&7&0&1&2&~\\
	\hline 5&4&~&~&1&0&~&~\\
	\hline 6&~&4&~&2&~&0&~\\
	\hline ~&~&~&~&~&~&~&~\\
	\hline 
	\end{tabular}
\end{tabular}
\end{center}

In general, we may take a latin square $L$ of order $n/2$ and form the
order $n$ latin square $L_1 \times L$  by defining:
\begin{align*}
L_1 \times L = &\{ (x,y;z), (x,y+n/2;z+n/2), (x+n/2,y;z+n/2),\\
	& (x+n/2,y+n/2;z) \mid (x,y;z) \in L \}
\end{align*}
The next Lemma is similar to the doubling
construction of 
\cite{MR84g:05036} 
which gives 2--critical sets.
\begin{lemma}
Let $M$ be a latin square of order $n$ such that 
$\mbox{gcs}(M)$ is $2$--critical. Then $\mbox{gcs}(L_1 \times M)$ 
is 2--critical.
\end{lemma}

\begin{proof}
Define $P$, a partial latin square of order $2n$, by
\begin{center}
	$P = ~$\begin{tabular}{|c|c|}
	\hline $M$ & $\mbox{gcs}^1(M)$ \\
	\hline $\mbox{gcs}^1(M)$ & $\mbox{gcs}(M)$ \\
	\hline
	\end{tabular}
\end{center}
Note that $P \subset L_1 \times M$.
Choose some $(i,j;k) \in \mathcal{Q}^n_{0,0}(P)$. Then $(i,j;k)$ will be
2--essential in one of two ways:

\begin{enumerate}

\item If $(i,j+n;k') \notin P$ then the set of cells 
\[
I= \{ (i,j;k),
(i,j+n;k'),
(i+n,j;k'),
(i+n,j+n;k) \}
\]
is an intercalate in $L_1 \times M$ such that 
$P \cap I = \{ (i,j;k) \}$.

\item Otherwise, $(i,j+n;k') \in P$. Since
$(i, j+n; k') \in \mbox{gcs}^1(M)$ which is
$2$--critical,
for some integers
$0 < \left| a \right|, \left| b \right| < n$
there exists an intercalate
\[
I = \{(i,j+n;k'), (i+a,j+n;l),(i,j+n+b;l),(i+a,j+n+b;k')\}
\]
for which
$I \cap \mbox{gcs}^1(M) = \{ (i,j+n;k') \}$. Hence there
is an intercalate
\[
I' = \{(i,j;k), (i,j+n+b;l), (i+n+a,j;l), (i+n+a,j+n+b;k)\}
\]
such that $I' \cap P = \{ (i,j;k) \}$ implying that
$(i,j;k)$ is 2--essential.
\end{enumerate}
Hence each $(i,j;k) \in P$ is the least element of an
intercalate so $P$ is 2--critical by Lemma~\ref{lemma:dmd-gcs}.
\end{proof}

\begin{corollary}
For all $s \geq 1$, $\mbox{gcs}(L_s) = P_s$ and $P_s$ is 2--critical.
\label{cor:P_s-2-critical}
\end{corollary}

\begin{lemma}
Let $\alpha$ be a row isotopism of $L_s$ defined below:
\begin{equation*}
\alpha(i) = 
	\begin{cases}
		4k_1+1, & i=4k_1+2 \\
		4k_1+2, & i=4k_1+1 \\
		\vdots &  \\
		4k_p+1, & i=4k_p+2 \\
		4k_p+2, & i=4k_p+1 \\
		i, & \text{otherwise}
	\end{cases}
\end{equation*}
where 
\begin{equation}
0 \leq p < 2^{s-2}\mbox{, }k_i \neq k_j \mbox{ for $i \neq j$, and $0
\leq 4k_i < 2^s$}
\label{eqn:innerswap}
\end{equation}
Then
$\mbox{gcs}(\alpha L_s)$ is 2--critical.
\label{lem:innerswap}
\end{lemma}

\begin{proof}
We proceed by induction. There are two base cases to check.
First, define
$H_2$ by the bracketed entries in the following square
and $\hat{H}_2$ to be the completion (as shown)
of $H_2$.

\begin{center}
	\begin{tabular}{|c|c|c|c|}
	\hline (0) & (1) & (2) & 3 \\
	\hline (2) & 3 & (0) & 1 \\
	\hline (1) & (0) & 3 & 2 \\
	\hline 3 & 2 & 1 & 0 \\
	\hline
	\end{tabular}
\end{center}
We note that $H_2$ is
isotopic to $P_2$ and so $H_2$ is a critical set. Further, each entry of
$H_2$ is the least element of some intercalate contained in
$\hat{H}_2$. For the second base case we need to
check a square of order $8$. First we construct
a general critical set $G_s$ of order $n = 2^s$ for $s \geq 3$ which will be shown to be
equivalent to $\mbox{gcs}(\alpha L_s)$ for $\alpha$
satisfying~\eqref{eqn:innerswap}.

Let $i$, $j$ be integers such that $i,j \equiv 0 \pmod{4}$ and $0 \leq i,j <
2^s$. We define each subsquare $\mathcal{Q}^4_{i,j}(G_s)$ as follows:
\begin{itemize}
\item If $\alpha(i+1) = i+1$ then set
$\mathcal{Q}^4_{i,j}(G_s) = \mathcal{Q}^4_{i,j}(P_s)$.

\item Otherwise, $\alpha(i+1) = i+2$,
$\alpha(i+2) = i+1$. 
Let $l$ be the integer such that
$\mathcal{Q}^4_{i,j}(L_s) = L_2^l$.
If 
$\mathcal{Q}^4_{i,j}(P_s)$ is similar to $L_2$
then set
$\mathcal{Q}^4_{i,j}(G_s) = \hat{H}^l_2$
otherwise
set $\mathcal{Q}^4_{i,j}(G_s) = H^l_2$ 
\end{itemize}

Since $H_2$ is isotopic to $P_2$ it follows that $G_s$
is isotopic to $P_s$, so $G_s$ is a critical set. 
To finish the second base case, we observe that
each entry of
$G_3$ is the least element of some intercalate contained in
$\alpha L_3$.

Next, fix the integer $s > 3$. 
We can partition $G_s$ into
$4 \times 4$ subsquares $\mathcal{Q}^4_{i,j}(G_s)$ where
$i,j \equiv 0 \pmod{4}$. 
There are two cases for each
subsquare:
\begin{enumerate}

\item If $\mathcal{Q}^4_{i,j}(G_s)$ is similar to $P_2$ or $H_2$
then each $(x,y;z) \in \mathcal{Q}^4_{i,j}(G_s)$ is
the least element of an intercalate.

\item Otherwise, $\mathcal{Q}^4_{i,j}(G_s)$ is isotopic to 
$\mathcal{Q}^4_{i,j}(P_s) = L_2$.  Since $G_s$ is isotopic to $P_s$, the
definition of $P_s$ implies that
$(i/4, j/4 ; k) \in P_{s-2}$ for some $k \in N_{n/4}$.
Then we know that there is an intercalate
\begin{align*}
I = \{
&(i/4, j/4; k),
(i/4 + a, j/4; k'), \\
&(i/4, j/4 + b; k'),
(i/4 + a, j/4 + b; k)
\}
\end{align*}
in $L_{s-2}$ such that
$I \cap P_{s-2} = \{ (i/4, j/4; k) \}$
and $a, b > 4$. Due to this intercalate $I$ and the 
the definition of $P_s$ we now see that
the subsquare
\[ R = \mathcal{Q}^4_{i,j}(G_s) \cup
\mathcal{Q}^4_{i,j+b+n/4}(G_s) \cup
\mathcal{Q}^4_{i+a+n/4,j}(G_s) \cup
\mathcal{Q}^4_{i+a+n/4,j+b+n/4}(G_s)
\]
is similar to $P_3$. We verified earlier that
each $(x,y;z) \in \mathcal{Q}^4_{i,j}(G_3)$ is the least element of an
intercalate and is 2--essential.
\end{enumerate}
By Lemma~\ref{lemma:dmd-gcs} we have
$G_s = \mbox{gcs}(\alpha L_s)$ and that $\mbox{gcs}(\alpha L_s)$ is 2--critical.
\end{proof}

\section{The Main Result}
\label{sec:main}
\begin{theorem}
Let $\alpha_{k,k'}$ be a row isotopism on a latin square of order $n = 2^s$,
defined by
\begin{equation*}
\alpha_{k,k'}(i) = 
	\begin{cases}
		k', & i=k \\
		k, & i=k' \\
		i, & \text{ otherwise.}
	\end{cases}
\end{equation*}
where 
\begin{equation}
\left| k-k' \right| < 3
\mbox{ and }
j \leq k < k' < j+4
\mbox{ for some }j \equiv 0 \pmod{4}
\label{eqn:allowableK}
\end{equation}
Then $\mbox{gcs}(\alpha_{k,k'}  L_s)$ is
2--critical, strong, and completes top down to $\alpha_{k,k'} L_s$.
\label{thm:cutdown}
\end{theorem}

The proof of Theorem~\ref{thm:cutdown} is based on induction. 
The case $s=2$ is treated separately in Section~\ref{sec:base-cases-L2}.
The remaining sections
contain the inductive proof,
beginning with the base case of $s=3$.

\subsection{Case $s=2$}
\label{sec:base-cases-L2}

The six possible $\mbox{gcs}(\alpha_{k,k'}L_2)$ are shown below.
Each critical set is
2--critical, strong, and completes top down.
\begin{align*}
\mbox{gcs}(\alpha_{0,1}L_2) = 
\begin{tabular}{|c|c|c|c|}
\hline  1  & (0) & (3) &  2  \\
\hline (0) & (1) &  2  &  3  \\
\hline (2) &  3  & (0) &  1  \\
\hline  3  &  2  &  1  &  0  \\
\hline
\end{tabular}
&
\quad \mbox{gcs}(\alpha_{0,2}L_2) = 
\begin{tabular}{|c|c|c|c|}
\hline (2) &  3  & (0) &  1  \\
\hline  1  & (0) & (3) &  2  \\
\hline (0) & (1) &  2  &  3  \\
\hline  3  &  2  &  1  &  0  \\
\hline
\end{tabular}\\
\mbox{gcs}(\alpha_{0,3}L_2) = 
\begin{tabular}{|c|c|c|c|}
\hline (3) & (2) & (1) &  0  \\
\hline (1) &  0  & (3) &  2  \\
\hline (2) & (3) &  0  &  1  \\
\hline  0  &  1  &  2  &  3  \\
\hline
\end{tabular}
&
\quad \mbox{gcs}(\alpha_{1,2}L_2) = 
\begin{tabular}{|c|c|c|c|}
\hline (0) & (1) & (2) &  3  \\
\hline (2) &  3  & (0) &  1  \\
\hline (1) & (0) &  3  &  2  \\
\hline  3  &  2  &  1  &  0  \\
\hline
\end{tabular} \\
\mbox{gcs}(\alpha_{1,3}L_2) = 
\begin{tabular}{|c|c|c|c|}
\hline (0) &  1  & (2) &  3  \\
\hline  3  & (2) & (1) &  0  \\
\hline (2) & (3) &  0  &  1  \\
\hline  1  &  0  &  3  &  2  \\
\hline
\end{tabular}
&
\quad \mbox{gcs}(\alpha_{2,3}L_2) = 
\begin{tabular}{|c|c|c|c|}
\hline  0  & (1) & (2) &  3  \\
\hline (1) & (0) &  3  &  2  \\
\hline (3) &  2  & (1) &  0  \\
\hline  2  &  3  &  0  &  1  \\
\hline
\end{tabular}
\end{align*}

\subsection{Base Cases for $s=3$}
\label{sec:base-cases-L3}

Let $s=3$ and $k, k' \in N_8$ such that~\eqref{eqn:allowableK}
is satisfied.
The base case for $s=3$ and $\alpha_{k,k'}$
is divided into two parts: $k,k' \geq 4$ and
$k,k' < 4$.
For each $(k,k')$ we see that the associated greedy critical set
$\mbox{gcs}(\alpha_{k,k'} L_3)$ is
2--critical, strong, and completes top down.

Let $\Gamma$ be the set of
$(k,k') \in N_8 \times N_8$ 
satisfying~\eqref{eqn:allowableK}
where $k,k' \geq 4$:
\begin{equation}
\Gamma = \{ (4,5), (4,6), (5,6), (5,7), (6,7) \}
\label{eqn:gamma}
\end{equation}
Suppose $(k,k') = (4,5)$. The
partial latin square
$\mbox{gcs}(\alpha_{4,5}L_3)$
is shown below as the entries in brackets. Also, we
take this opportunity to define the partial latin square
$E(4,5)_2$.
\begin{align*}
\mbox{gcs}(\alpha_{4,5}L_3)
& = \begin{tabular}{|c|c|c|c||c|c|c|c|}
	\hline (0) & (1) &  2  & (3) & (4) & (5) & (6) &  7  \\
	\hline (1) & (0) & (3) & (2) & (5) & (4) &  7  &  6  \\
	\hline (2) & (3) & (0) & (1) & (6) &  7  & (4) &  5  \\
	\hline (3) & (2) & (1) & (0) &  7  &  6  &  5  &  4  \\
	\hline
	\hline  5  & (4) & (7) &  6  &  1  & (0) & (3) &  2  \\
	\hline (4) & (5) &  6  &  7  & (0) & (1) &  2  &  3  \\
	\hline (6) &  7  & (4) &  5  & (2) &  3  & (0) &  1  \\
	\hline  7  &  6  &  5  &  4  &  3  &  2  &  1  &  0  \\
	\hline
	\end{tabular}
=
	\begin{tabular}{|c||c|}
	\hline $E(4,5)_2$ & $\bullet$ \\
	\hline
	\hline $\bullet$ & $\bullet$ \\
	\hline
	\end{tabular}
\end{align*}
The other 
$\mbox{gcs}(\alpha_{k,k'}L_3)$
and
$E(k,k')_2$ 
for
$(k,k') \in \Gamma$
are shown in Appendix~\ref{app:A}.

Otherwise, $k,k' < 4$.
Let
$\Lambda = \{ (0,1), (0,2),
(0,3),
(1,2),
(1,3),
(2,3) \}$.
Each case defines
a partial latin square
$A(k,k')_2$. For example,
\begin{align*}
\mbox{gcs}(\alpha_{0,1}L_3) 
&= 	\begin{tabular}{|c|c|c|c||c|c|c|c|}
	\hline  1  & (0) & (3) & (2) &  5  & (4) & (7) &  6  \\
	\hline (0) & (1) & (2) & (3) & (4) & (5) &  6  &  7  \\
	\hline (2) & (3) & (0) & (1) & (6) &  7  & (4) &  5  \\
	\hline (3) & (2) & (1) & (0) &  7  &  6  &  5  &  4  \\
	\hline
	\hline (4) & (5) & (6) &  7  & (0) & (1) & (2) &  3  \\
	\hline (5) & (4) &  7  &  6  & (1) & (0) &  3  &  2  \\
	\hline (6) &  7  & (4) &  5  & (2) &  3  & (0) &  1  \\
	\hline  7  &  6  &  5  &  4  &  3  &  2  &  1  &  0  \\
	\hline
	\end{tabular}
= \begin{tabular}{|c||c|}
	\hline $A(0,1)_2$ & $\bullet$ \\
	\hline
	\hline $\bullet$ & $\bullet$ \\
	\hline
	\end{tabular}
\end{align*}
The other
$A(k,k')_2$
are shown in Appendix~\ref{app:B}.
\subsection{The Final Construction}

In this subsection we will define partial latin squares
$E(k,k')_{s}$,
$A(k,k')_{s}$,
and
$G(k,k', s)$.  The squares
$E(k,k')_{s}$
and
$A(k,k')_{s}$
are used in recursively defining
$G(k,k', s)$.
In Section~\ref{sec:completing-the-proof} we will show
that 
$G(k,k', s) = \mbox{gcs}(\alpha_{k,k'} L_s)$.

Recall that
$A(k,k')_{2}$
and
$E(k,k')_{2}$ 
were defined in Section~\ref{sec:base-cases-L3}.
For $s \geq 4$ 
and
$\delta=2^{s-2}$
if $k,k' < s^{s-1}$
we define
\begin{equation}
A(k,k')_{s-1} = 
	\begin{cases}
		\begin{tabular}{|c|c|}
		\hline $A(k,k')^0_{s-2}$ & $A(k,k')^1_{s-2}$ \\
		\hline $L^1_{s-2}$ & $L^0_{s-2}$ \\
		\hline
		\end{tabular} &\mbox{ if } k,k' \in N_{\delta}^0 \\

		~ \\

		\begin{tabular}{|c|c|}
		\hline $L^0_{s-2}$ & $L^1_{s-2}$ \\
		\hline $A(k-\delta,k'-\delta)^1_{s-2}$ & $A(k-\delta,k'-\delta)^0_{s-2}$ \\
		\hline
		\end{tabular} &\mbox{ if } k,k' \in N_{\delta}^1
	\end{cases}
\end{equation}
\begin{equation}
E(k,k')_{s-1} = 
	\begin{cases}
		\begin{tabular}{|c|c|}
		\hline $L^0_{s-2}$ & $E(k-\delta,k'-\delta)^1_{s-2}$ \\
		\hline $L^1_{s-2}$ & $L^0_{s-2}$ \\
		\hline
		\end{tabular} &\mbox{ if } k,k' \in N_{\delta}^2 \\

		~ \\

		\begin{tabular}{|c|c|}
		\hline $E(k-2\delta,k'-2\delta)^0_{s-2}$ & $L^1_{s-2}$ \\
		\hline $E(k-2\delta,k'-2\delta)^1_{s-2}$ & $E(k-2\delta,k'-2\delta)^0_{s-2}$ \\
		\hline
		\end{tabular} &\mbox{ if } k,k' \in N_{\delta}^3
	\end{cases}
\label{eqn:E}
\end{equation}
The previous section, Appendix~\ref{app:A}, and Appendix~\ref{app:B}
give $G(k,k', s) = \mbox{gcs}(\alpha_{k,k'} L_3)$. For $s \geq 4$,
define
\begin{equation}
G(k,k', s) = \\
	\begin{cases}
		\begin{tabular}{|c|c|}
		\hline $A(k,k')^0_{s-1}$ & $G(k,k', s-1)$ \\
		\hline $P^0_{s-1}$ & $P^1_{s-1}$ \\
		\hline
		\end{tabular} &\mbox{ if } k, k' < n/2 \\

		~ \\

		\begin{tabular}{|c|c|}
		\hline $E(k,k')^0_{s-1}$ & $P^1_{s-1}$ \\
		\hline $G(k-2\delta,k'-2\delta, s-1)^1$ & $G(k-2\delta,k'-2\delta, s-1)^0$ \\
		\hline
		\end{tabular} &\mbox{ if } k,k' \geq n/2
	\end{cases}
\label{eqn:Gkk}
\end{equation}

The following Lemma is immediate from the
definition of $E(k,k')_2$ and \eqref{eqn:E}.

\begin{lemma}
Let $i,j \equiv 0 \pmod{4}$ and
$W = \mathcal{Q}^4_{i,j}(E(k,k')_s)$. Then $W \approx L_2$
or $W \approx E(l,l')_2$ for some
$(l,l') \in \Gamma$.
\label{lem:Estructure}
\end{lemma}

For each $(k,k') \in \Gamma$, define
\begin{equation}
	U(k,k') = \begin{tabular}{|c|c|}
	\hline $E(k,k')^0_2$ & $P^1_2$ \\
	\hline $E(k,k')^1_2$ & $P^0_2$\\
	\hline
	\end{tabular}
\end{equation}

\begin{lemma}
Let $r,c \in N_4$. Then 
for each $(k,k') \in \Gamma$,
\[ \mathcal{A}(U(k,k'),N_8 \times N_8,L_3)_{r,c} \cap N^1_4 = \emptyset \]
\label{lem:blocking1}
\end{lemma}

\begin{proof}
There are four cases to inspect. 
First, let $U=U(4,5) =U(4,6)$.
\begin{equation*}
	U = 
	\begin{tabular}{|c|c|c|c||c|c|c|c|}
	\hline (0) & (1) &  2  & (3) & (4) & (5) & (6) &  7  \\
	\hline (1) & (0) & (3) & (2) & (5) & (4) &  7  &  6  \\
	\hline (2) & (3) & (0) & (1) & (6) &  7  & (4) &  5  \\
	\hline (3) & (2) & (1) & (0) &  7  &  6  &  5  &  4  \\
	\hline
	\hline (4) & (5) &  6  & (7) & (0) & (1) & (2) &  3  \\
	\hline (5) & (4) & (7) & (6) & (1) & (0) &  3  &  2  \\
	\hline (6) & (7) & (4) & (5) & (2) &  3  & (0) &  1  \\
	\hline (7) & (6) & (5) & (4) &  3  &  2  &  1  &  0  \\
	\hline 
	\end{tabular}
\end{equation*}
In the first row there are 
empty cells $(0,2)$ and $(0,7)$ which could be filled with a $2$ or $7$. 
However $7 \notin \mathcal{A}(U,N_8 \times N_8,L_3)_{0,2}$ 
since
$(5,2;7) \in U$.
So
$\mathcal{A}(U,N_8 \times N_8,L_3)_{0,2} = \{2\} \subset N_4^0$. Finally,
$\mathcal{A}(U,N_8 \times N_8,L_3)_{r,c} = \emptyset$ for $r,c \in N_4^0$ 
and $(r,c) \neq (0,2)$.
This completes the proof of this case.
The other three $U(k,k')$ are displayed in Appendix~\ref{app:0}.
\end{proof}

For each $(k,k') \in \Gamma$, define
\begin{equation}	
	V(k,k') = \begin{tabular}{|c|c|}
	\hline $E(k,k')^0_2$ & $P^1_2$ \\
	\hline $P^1_2$ & $E(k,k')^0_2$  \\
	\hline
	\end{tabular}
\end{equation}	

\begin{lemma}
Let
$r \in N_4$, 
$c \in N^1_4$. Then for each
$(k,k') \in \Gamma$
\[ \mathcal{A}(V(k,k'),N_8 \times N_8,L_3)_{r,c} \cap N^0_4 = \emptyset \]
\label{lem:blocking2}
\end{lemma}
As with the previous lemma there are four cases to check with very
similar reasoning (see Appendix~\ref{app:1}).

\subsection{Completing the Proof of Theorem~\ref{thm:cutdown}}
\label{sec:completing-the-proof}

The proof of Theorem~\ref{thm:cutdown} will require a few technical
Lemmas. First, Lemma~\ref{lem:intercalatesEverywhere} follows directly from the definition
of $L_s$.

\begin{lemma}
Let $(i,j;k)$, $(i,j';k') \in L_s$ where
$n = 2^s$ and
\begin{align*}
&0 \leq i < n/2 \\
&0 \leq j < n/2 \leq j' < n
\end{align*}
Then there exists an integer $i'$ with $n/2 \leq i' < n$ such that 
\[
\{ (i,j;k),	(i,j';k'),
(i',j;k'),	(i',j';k) \}
\]
is an intercalate
in $L_s$.
\label{lem:intercalatesEverywhere}
\end{lemma}
Using ideas in the proof of Lemma~\ref{lem:innerswap}, we also have:
\begin{corollary}
Let $(i,j;k)$, $(i,j';k') \in L_s$ where
$n = 2^s$, $i,j,j' \equiv 0 \pmod{4}$ and
\begin{align*}
&0 \leq i < n/2 \\
&0 \leq j < n/2 \leq j' < n
\end{align*}
Then there exists an integer $i' \equiv \pmod{4}$ with $n/2 \leq i' < n$
such that
\[
\mathcal{Q}^4_{i,j}(L_s) \cup 
\mathcal{Q}^4_{i,j'}(L_s) \cup 
\mathcal{Q}^4_{i',j}(L_s) \cup 
\mathcal{Q}^4_{i',j'}(L_s)
\approx L_3
\]
\label{cor:intercalatesEverywhere}
\end{corollary}
\begin{lemma}
Let $P$ be a partial latin square contained in the latin square $L$ of
order $n$.  Let $R_1, R_2, C \subseteq N_n$ and define 
$S = (R_1 \cup R_2) \times C$
so that
\begin{align*}
\left| C \right| &= \left| R_1 \cup R_2 \right| \\
Q &= \{ (i,j;k) \mid (i,j) \in S \mbox{, } (i,j;k) \in P \} \\
L' &= \{ (i,j;k) \mid (i,j) \in S \mbox{, } (i,j;k) \in L \}
\end{align*}
where $L'$ is a latin subsquare of $L$.
If $Q$ strongly completes (extends) to 
\[
\hat{Q} = Q \cup \{ (i,j;k) \mid (i,j;k) \in L' \mbox{, $i \in R_1$, $c \in C$} \}
\]
and
$\mathcal{A}(P, R_1 \times C, L) = \mathcal{A}(Q, R_1 \times C, L')$ 
then $P$ has a unique extension to 
\begin{equation}
\hat{P} = P \cup \{ (i,j;k) \mid (i,j;k) \in L \mbox{, $i \in R_1$, $c \in C$} \}.
\label{eqn:P-completion}
\end{equation}
\label{lem:simpler-subsquare-uc}
\end{lemma}

\begin{remark}
The subsquare $Q$ has strong completion through rows $R_1$ only. This is
useful if rows $R_2$ of the arrays of alternatives are not equal. On the
other hand, if we set $R_2 = \emptyset$ then the lemma says that
$P$ extends to $P \cup L'$.
\end{remark}

\begin{proof}[Proof of Lemma~\ref{lem:simpler-subsquare-uc}]
Since $Q$ has strong completion through rows $R_1$ in $L'$ there must be a sequence
\[
(Q^{(1)}, r_1, c_1),
(Q^{(2)}, r_2, c_2),
\ldots,
(Q^{(m)}, r_m, c_m)
\]
such that
$Q^{(1)} = Q$,
$Q^{(m+1)} = \hat{Q}$,
and
$Q^{(i)} \subsetneq Q^{(i+1)}$,
$r_i \in R_1$,
$c_i \in C$
for 
each $i$.
Also, the pairs $(r_i,c_i)$ are distinct and 
$\left| \mathcal{A}(Q^{(i)}, S,L')_{{r_i}, {c_i}} \right| = 1$ for each $i$.

Define the sequence $(P^{(i)}, r_i, c_i) = (P \cup Q^{(i)}, r_i, c_i)$ for $1
\leq i \leq m$. It is obvious that $P^{(i)} \subsetneq P^{(i+1)}$. To show
that this is a strong completion we need 
$\left| \mathcal{A}(P^{(i)}, S,L)_{{r_i}, {c_i}} \right| = 1$ for each $i$.
First, $\left| \mathcal{A}(P^{(1)}, S,L)_{{r_1}, {c_1}} \right| = 1$ 
by the definition of $Q$.
Now suppose that
$\left| \mathcal{A}(P^{(i)}, S,L)_{{r_{i}}, {c_{i}}} \right| = 1$ and
\[
\mathcal{A}(P^{(i)}, S, L)_{r,c} = \mathcal{A}(Q^{(i)}, S, L')_{r,c} \mbox{ for each $r \in
R_1$, $c \in C$}
\]
where $i > 1$ is fixed.
Fill the cell $(r_i,c_i)$ in $Q^{(i)}$ with the (unique) symbol 
$e_i \in \mathcal{A}(Q^{(i)}, S, L')_{{r_i},{c_i}}$ to get $Q^{(i+1)}$. Do the same in
$P^{(i)}$. Now let
\begin{equation*}
\mathcal{A}(P^{(i+1)}, S, L)_{r,c} = 
	\begin{cases}
		\mathcal{A}(P^{(i)}, S, L)_{r,c}, &\mbox{ if $r \neq r_i$ and $c \neq c_i$} \\
		\mathcal{A}(P^{(i)}, S, L)_{r,c} \setminus \{ e_i \}, &\mbox{ otherwise}
	\end{cases}
\end{equation*}
Since $\mathcal{A}(Q^{(i+1)}, S, L)_{r,c}$ has the same definition (i.e. the
symbol $e_i$ deleted from the corresponding row and column) it follows that
$\mathcal{A}(P^{(i+1)}, S, L)_{r,c} = \mathcal{A}(Q^{(i+1)}, S, L)_{r,c}$
for each $r$, $c$. Hence
$P$ has strong completion to $\hat{P}$.
\end{proof}

\begin{lemma}
$G(k,k',s)$ has strong completion top down.
\label{lem:G-strong-top-down}
\end{lemma}

\begin{proof}
The case $s=2$ and base case $s=3$ are given earlier. So suppose that the theorem is true for all 
$L_t$ where $3 \leq t < s$.
Let $\delta = 2^{s-2}$. 
There are four cases depending on where the row
swap occurs.

Case~1: 
$k,k' \in N_{\delta}^3$. Write 
$l = k-2\delta$,
$l' = k'-2\delta$,
$h = k-3\delta$,
$h' = k'-3\delta$.
Then
\begin{align}
&G(k,k', s) =
		\begin{tabular}{|c|c|}
		\hline $E(k,k')^0_{s-1}$ & $P^1_{s-1}$ \\
		\hline $G(l,l', s-1)^1$ & $G(l,l', s-1)^0$ \\
		\hline
		\end{tabular} \nonumber \\
	&=
	\begin{tabular}{|c|c||c|c|}
	\hline $E(l,l')^0_{s-2}$ & $L^1_{s-2}$ & $L^2_{s-2}$ & $P^3_{s-2}$ \\
	\hline $E(l,l')^1_{s-2}$ & $E(l,l')^0_{s-2}$ & $P^3_{s-2}$ & $P^2_{s-2}$ \\
	\hline
	\hline $E(l,l')^2_{s-2}$ & $P^3_{s-2}$ & $E(l,l')^0_{s-2}$ & $P^1_{s-2}$ \\
	\hline $G(h,h', s-2)^3$ & $G(h,h', s-2)^2$ & $G(h,h', s-2)^1$ & $G(h,h', s-2)^0$ \\
	\hline
	\end{tabular}
\label{eq:case1}
\end{align}
First we show that the top $\delta$ rows of
$G(k,k', s)$ 
have unique completion top down. Since the cells
$N_{\delta}^0 \times (N_{\delta}^1 \cup N_{\delta}^2)$
are completely filled in, we need to show that the cells
$N_{\delta}^0 \times (N_{\delta}^0 \cup N_{\delta}^3)$
have unique completion. We will use
Lemma~\ref{lem:simpler-subsquare-uc}. Let
\begin{align*}
R_1 &= N_{\delta}^0, \quad R_2 = N_{\delta}^3, \quad C = N_{\delta}^0
\cup N_{\delta}^3,
\quad L = \alpha_{k,k'} L_s, \quad P = G(k,k',s)
\end{align*}
Then $Q$ and $L'$ are defined to be
\begin{align*}
Q &= 	\begin{tabular}{|c|c|}
	\hline $E(l,l')^0_{s-2}$ & $P^3_{s-2}$ \\
	\hline $G(h,h', s-2)^3$ & $G(h,h', s-2)^0$ \\
	\hline
	\end{tabular} 
= G(l,l', s-1)  \\
L' &= 	\begin{tabular}{|c|c|}
	\hline $L^0_{s-2}$ & $L^3_{s-2}$ \\
	\hline $L^3_{s-2}$ & $L^0_{s-2}$ \\
	\hline
	\end{tabular}
= \alpha_{l,l'} L_{s-1}
\end{align*}
By the inductive hypothesis $Q$ has
unique completion top down to $L'$.  
From~\eqref{eq:case1} we see that
$\mathcal{A}(P, R_1 \times C, L) \cap (N_{\delta}^1 \cup N_{\delta}^2) = \emptyset$
and
$\mathcal{A}(P, R_1 \times C, L) = \mathcal{A}(Q, R_1 \times C, L')$.
Now Lemma~\ref{lem:simpler-subsquare-uc} gives the strong top down 
extension
$G^+(k,k',s)$ of $G(k,k',s)$:
\begin{align*}
	\begin{tabular}{|c|c||c|c|}
	\hline $L^0_{s-2}$ & $L^1_{s-2}$ & $L^2_{s-2}$ & $L^3_{s-2}$ \\
	\hline $E(l,l')^1_{s-2}$ & $E(l,l')^0_{s-2}$ & $P^3_{s-2}$ & $P^2_{s-2}$ \\
	\hline
	\hline $E(l,l')^2_{s-2}$ & $P^3_{s-2}$ & $E(l,l')^0_{s-2}$ & $P^1_{s-2}$ \\
	\hline $G(h,h', s-2)^3$ & $G(h,h', s-2)^2$ & $G(h,h', s-2)^1$ & $G(h,h', s-2)^0$ \\
	\hline
	\end{tabular}
\end{align*}
We will now show that
\begin{align}
&\mathcal{A}(G^{+}(k,k',s), N_n \times N_n, \alpha_{k,k'} L_s)_{r,c} \subseteq 
N^1_{\delta} \cup
N^3_{\delta}
\mbox{ for $r \in N_{\delta}^1$, $c \in N_{\delta}^0$}  \label{eqn:restrictingAA1} \\
&\mathcal{A}(G^{+}(k,k',s), N_n \times N_n, \alpha_{k,k'} L_s)_{r,c} \subseteq 
N^1_{\delta} \cup
N^3_{\delta}
\mbox{ for $r \in N_{\delta}^1$, $c \in N_{\delta}^2$}
\label{eqn:restrictingAA2}
\end{align}
Let $i$, $j$, $u$, $v$ be integers such that $i,j,u,v \equiv 0
\pmod{4}$ and
\begin{align*}
i, i+3 &\in N_{\delta}^1 \quad \quad
j, j+3 \in N_{\delta}^0 \quad \quad v, v+3 \in N_{\delta}^3
\end{align*}
Then by Lemma~\ref{lem:Estructure},
$\mathcal{Q}^4_{i,j}(G^{+}(k,k',s)) \supseteq E(k,k')^w_2$ for
some $(k,k')$ and integer $w$.
Also,
$\mathcal{Q}^4_{u,v}(G^{+}(k,k',s)) \supseteq P_2^x$ for some integer $x$. So to
restrict the array of alternatives we apply
Lemma~\ref{lem:blocking1}
to the subsquare (which exists due to
Corollary~\ref{cor:intercalatesEverywhere}
for some $x \in N_{\delta}^2$)
\[
\mathcal{Q}^4_{i,j}(G^{+}(k,k',s))
\cup
\mathcal{Q}^4_{i,v}(G^{+}(k,k',s))
\cup
\mathcal{Q}^4_{x,j}(G^{+}(k,k',s))
\cup
\mathcal{Q}^4_{x,v}(G^{+}(k,k',s))
\]
and
Lemma~\ref{lem:blocking2}
to the subsquare (which exists due to
Corollary~\ref{cor:intercalatesEverywhere}
for some $y \in N_{\delta}^2$)
\[
\mathcal{Q}^4_{i,j + \delta}(G^{+}(k,k',s))
\cup
\mathcal{Q}^4_{i,v-\delta}(G^{+}(k,k',s))
\cup
\mathcal{Q}^4_{y,j + \delta}(G^{+}(k,k',s))
\cup
\mathcal{Q}^4_{y,v-\delta}(G^{+}(k,k',s))
\]
which gives
\begin{align*}
\mathcal{A}(G^{+}(k,k',s), N_n \times N_n, \alpha_{k,k'} L_s)_{r,c} &\cap 
N_{\delta}^2
= \emptyset \quad  \mbox{ for $r \in N_{\delta}^1$, $c \in N_{\delta}^0$} \\
\mathcal{A}(G^{+}(k,k',s), N_n \times N_n, \alpha_{k,k'} L_s)_{r,c} &\cap N_{\delta}^0
= \emptyset \quad \mbox{ for $r \in N_{\delta}^1$, $c \in N_{\delta}^2$}
\end{align*}
which imply~\eqref{eqn:restrictingAA1} and \eqref{eqn:restrictingAA2}.
Now consider these two subsquares of $G^{+}(k,k', s)$:
\begin{align}
	&\begin{tabular}{|c|c||c|c|}
	\hline \phantom{$L^0_{s-2}$} & \phantom{$L^1_{s-2}$} & \phantom{$L^2_{s-2}$} & \phantom{$L^3_{s-2}$} \\
	\hline $E(l,l')^1_{s-2}$ &             & $P^3_{s-2}$ &             \\
	\hline
	\hline             &             &             &             \\
	\hline $G(h,h',{s-2})^3$ &             & $G(h,h',{s-2})^1$ &             \\
	\hline
	\end{tabular} \label{eq:left} \\
~ \nonumber \\
	&\begin{tabular}{|c|c||c|c|}
	\hline \phantom{$L^0_{s-2}$} & \phantom{$L^1_{s-2}$} & \phantom{$L^2_{s-2}$} & \phantom{$L^3_{s-2}$} \\
	\hline             & $E(l,l')^0_{s-2}$ &             & $P^2_{s-2}$ \\
	\hline
	\hline             &             &             &             \\
	\hline             & $G(h,h',{s-2})^2$ &          &  $G(h,h',{s-2})^0$  \\
	\hline
	\end{tabular} \label{eq:right}
\end{align}

We can now interleave the application of Lemma~\ref{lem:simpler-subsquare-uc}
with the inductive hypothesis
to show that the top halves of these subsquares completes strongly top
down. Suppose row $r$, for $r \in N_{\delta}^1$ has been completed 
in~\eqref{eq:left}. Then the array of
alternatives for subsquares identified in~\eqref{eq:right} is restricted such that the
inductive hypothesis applies. In other words,
\begin{align*}
&\mathcal{A}(G^{+}(k,k',s), N_n \times N_n, \alpha_{k,k'} L_s)_{r,c} \subseteq 
N^0_{\delta} \cup
N^2_{\delta}
\mbox{ for $c \in N_{\delta}^1$} \\
&\mathcal{A}(G^{+}(k,k',s), N_n \times N_n, \alpha_{k,k'} L_s)_{r,c} \subseteq 
N^0_{\delta} \cup
N^2_{\delta}
\mbox{ for $c \in N_{\delta}^3$}
\end{align*}
Hence $G^{+}(k,k',s)$ strongly extends top down to
$G^{++}(k,k',s)$:
\begin{align*}
	\begin{tabular}{|c|c||c|c|}
	\hline $L^0_{s-2}$ & $L^1_{s-2}$ & $L^2_{s-2}$ & $L^3_{s-2}$ \\
	\hline $L^1_{s-2}$ & $L^0_{s-2}$ & $L^3_{s-2}$ & $L^2_{s-2}$ \\
	\hline
	\hline $E(l,l')^2_{s-2}$ & $P^3_{s-2}$ & $E(l,l')^0_{s-2}$ & $P^1_{s-2}$ \\
	\hline $G(h,h', s-2)^3$ & $G(h,h', s-2)^2$ & $G(h,h', s-2)^1$ & $G(h,h', s-2)^0$ \\
	\hline
	\end{tabular}
\end{align*}
Finally, interleave the application of Lemma~\ref{lem:simpler-subsquare-uc}
with $R_2 = \emptyset$
to the subsquares
\[
R_1 \times C = \{ 2\delta, 2\delta+1, \ldots, 4\delta-1 \} \times \{ 0, 1, \ldots, 2\delta-1\}
\]
and
\[
R_1 \times C = \{ 2\delta, 2\delta+1, \ldots, 4\delta-1 \} \times \{ 2\delta, 2\delta+1, \ldots, 4\delta-1 \}
\]
which finishes the completion of $G(k,k',s)$.

Case~2: 
$k,k' \in N_{\delta}^2$.
Write 
$l = k-\delta$,
$l' = k'-\delta$,
$h = k - 2\delta$,
$h' = k' - 2\delta$. Then
\begin{align*}
	G(k,k',s) &= \begin{tabular}{|c|c||c|c|}
	\hline $L^0_{s-2}$ & $E(l,l')^1_{s-2}$ & $L^2_{s-2}$ & $P^3_{s-2}$ \\
	\hline $L^1_{s-2}$ & $L^0_{s-2}$ & $P^3_{s-2}$ & $P^2_{s-2}$ \\
	\hline
	\hline $A(h,h')^2_{s-2}$ & $G(h, h',s-2)^3$ & $A(h,h')^0_{s-2}$ & $G(h, h',s-2)^1$ \\
	\hline $P^3_{s-2}$ & $P^2_{s-2}$ & $P^1_{s-2}$ & $P^0_{s-2}$ \\
	\hline
	\end{tabular}
\end{align*}
and the reasoning is simpler than Case~1.

Case~3: 
$k,k' \in N_{\delta}^1$. Write
$l = k-\delta$,
$l' = k'-\delta$.

\begin{align*}
	G(k,k',s) = \begin{tabular}{|c|c||c|c|}
	\hline $L^0_{s-2}$ & $L^1_{s-2}$ & $E(k,k')^2_{s-2}$ & $P^3_{s-2}$ \\
	\hline $A(l,l')^1_{s-2}$ & $A(l,l')^0_{s-2}$ & $G(l,l',s-2)^3$ & $G(l,l',s-2)^2$ \\
	\hline
	\hline $L^2_{s-2}$ & $P^3_{s-2}$ & $L^0_{s-2}$ & $P^1_{s-2}$ \\
	\hline $P^3_{s-2}$ & $P^2_{s-2}$ & $P^1_{s-2}$ & $P^0_{s-2}$ \\
	\hline
	\end{tabular}
\end{align*}

Case~4: 
$k,k' \in N_{\delta}^0$.
\begin{align*}
	G(k,k',s) = \begin{tabular}{|c|c||c|c|}
	\hline $A(k,k')^0_{s-2}$ & $A(k,k')^1_{s-2}$ & $A(k,k')^2_{s-2}$ & $G(k,k',s-2)^3$ \\
	\hline $L^1_{s-2}$ & $L^0_{s-2}$ & $P^3_{s-2}$ & $P^2_{s-2}$ \\
	\hline
	\hline $L^2_{s-2}$ & $P^3_{s-2}$ & $L^0_{s-2}$ & $P^1_{s-2}$ \\
	\hline $P^3_{s-2}$ & $P^2_{s-2}$ & $P^1_{s-2}$ & $P^0_{s-2}$ \\
	\hline
	\end{tabular}
\end{align*}
This completes the proof.
\end{proof}

\begin{lemma}
$G(k,k',s)$ is a 2--critical greedy critical set.
\label{lem:2crit-gcs}
\end{lemma}

\begin{proof}
The case $s=2$ and base case $s=3$ were given earlier.
Suppose that the result is true for all 
$G(k,k',t)$ where $3 \leq t < s$. Write $\delta = 2^{s-2}$. There are
four cases for $k,k'$.
First, suppose that
$k,k' \in N_{\delta}^0$. Write down $G = G(k,k',s)$:
\begin{align*}
G &= 
	\begin{tabular}{|c|c||c|c|}
	\hline $A(k,k')^0_{s-2}$ & $A(k,k')^1_{s-2}$ & $A(k,k')^2_{s-2}$ & $G(k,k',s-2)^3$ \\
	\hline $L^1_{s-2}$ & $L^0_{s-2}$ & $P^3_{s-2}$ & $P^2_{s-2}$ \\
	\hline
	\hline $L^2_{s-2}$ & $P^3_{s-2}$ & $L^0_{s-2}$ & $P^1_{s-2}$ \\
	\hline $P^3_{s-2}$ & $P^2_{s-2}$ & $P^1_{s-2}$ & $P^0_{s-2}$ \\
	\hline
	\end{tabular}
\end{align*}
Now identify the following subsquares:
\begin{enumerate}
	\item $\mathcal{Q}^{2\delta}_{2\delta, 0}(G) \approx P_{s-1}$
	
	\item $\mathcal{Q}^{2\delta}_{2\delta, 2\delta}(G) \approx P_{s-1}$
	
	\item $\mathcal{Q}^{2\delta}_{0, 2\delta}(G) \approx G(k,k',s-1)$

	\item $	\mathcal{Q}^{\delta}_{0, 0}(G) \cup
		\mathcal{Q}^{\delta}_{0, 3\delta}(G) \cup
		\mathcal{Q}^{\delta}_{3\delta, 0}(G) \cup
		\mathcal{Q}^{\delta}_{3\delta, 3\delta}(G)
		\approx G(k,k',s-1)$

	\item $	\mathcal{Q}^{\delta}_{0, \delta}(G) \cup
		\mathcal{Q}^{\delta}_{0, 3\delta}(G) \cup
		\mathcal{Q}^{\delta}_{2\delta, \delta}(G) \cup
		\mathcal{Q}^{\delta}_{2\delta, 3\delta}(G)
		\approx G(k,k',s-1)$
	
	\item $	\mathcal{Q}^{\delta}_{\delta, 0}(G) \cup
		\mathcal{Q}^{\delta}_{\delta, 2\delta}(G) \cup
		\mathcal{Q}^{\delta}_{3\delta, 0}(G) \cup
		\mathcal{Q}^{\delta}_{3\delta, 2\delta}(G)
		\approx P_{s-1}$
	
	\item $	\mathcal{Q}^{\delta}_{\delta, \delta}(G) \cup
		\mathcal{Q}^{\delta}_{\delta, 3\delta}(G) \cup
		\mathcal{Q}^{\delta}_{3\delta, \delta}(G) \cup
		\mathcal{Q}^{\delta}_{3\delta, 3\delta}(G)
		\approx P_{s-1}$

\end{enumerate}
With these subsquares,
the inductive hypothesis, and Corollary~\ref{cor:P_s-2-critical}, 
we see that each entry $x \in G(k,k',s)$ is
2--essential and that there exists a trade $T_{x} \subset
\alpha_{k,k'}
L_s$ such that $G(k,k',s) \cap T_{x} = \{ x \}$ and
$x$ is the least element of $T_{x}$. 
Further, $G(k,k',s)$ has strong completion 
by Lemma~\ref{lem:G-strong-top-down} so 
$G(k,k',s)$ is a critical set. Finally,
Lemma~\ref{lemma:dmd-gcs} shows that
$G = \mbox{gcs}(\alpha_{k,k'} L_s)$.
We omit the remaining three cases where $k,k'$ are in  $N_{\delta}^1, N_{\delta}^2, N_{\delta}^3$ since
the reasoning is very similar.
\end{proof}

\begin{theoremMain}
Let $\alpha_{k,k'}$ be a row isotopism on a latin square of order $n = 2^s$,
defined by
\begin{equation*}
\alpha_{k,k'}(i) = 
	\begin{cases}
		k', & i=k \\
		k, & i=k' \\
		i, & \text{ otherwise.}
	\end{cases}
\end{equation*}
where 
\begin{equation}
\left| k-k' \right| < 3
\mbox{ and }
j \leq k < k' < j+4
\mbox{ for some }j \equiv 0 \pmod{4}
\end{equation}
Then $\mbox{gcs}(\alpha_{k,k'}  L_s)$ is
2--critical, strong, and completes top down to $\alpha_{k,k'} L_s$.
\end{theoremMain}

\begin{proof}
The Theorem
follows from 
Lemmas~\ref{lem:G-strong-top-down}
and 
\ref{lem:2crit-gcs}.
\end{proof}

\section{Conclusion}

We believe that a stronger version of the
theorem is true, where~\eqref{eqn:allowableK} is weakened.

\begin{conjecture}
Let $\alpha_{k,k'}$ be a row isotopism on a latin square of order $n = 2^s$,
defined by
\begin{equation*}
\alpha(i) = 
	\begin{cases}
		k', & i=k \\
		k, & i=k' \\
		i, & \text{ otherwise.}
	\end{cases}
\end{equation*}
where $\left| k-k' \right| < 3$
Then $\mbox{gcs}(\alpha_{k,k'}  L_s)$ is
2--critical, strong, and completes top down.
\end{conjecture}

We have verified the conjecture by computer search for $2 \leq s \leq 5$
and all possible $\alpha_{k,k'}$.

\clearpage

\appendix
\makeatletter
\def\@seccntformat#1{\csname Pref@#1\endcsname \csname the#1\endcsname\quad}
\def\Pref@section{Appendix~}
\makeatother

\section{~}
\label{app:G60example}

Suppose we wish to calculate $\mbox{gcs}(\alpha_{60,62}L_{6})$. 

\begin{itemize}
\item First, $s=6$ so
$n=2^5 = 64$, and $\delta = 2^{s-2} = 2^4 = 16$.
Use \eqref{eqn:Gkk} to write down $G(60,62,6)$:
\begin{align*}
G(60,62,6) 
&= \begin{tabular}{|c|c|}
	\hline $E(60,62)^0_{5}$ & $P^1_{5}$ \\
	\hline $G(60-32,62-32, 5)^1$ & $G(60-32,62-32, 5)^0$ \\
	\hline
	\end{tabular} \\
&= \begin{tabular}{|c|c|}
	\hline $E(60,62)^0_{5}$ & $P^1_{5}$ \\
	\hline $G(28,30, 5)^1$ & $G(28,30, 5)^0$ \\
	\hline
	\end{tabular}
\end{align*}

\item Now use \eqref{eqn:Gkk} twice more:
\begin{align*}
G(28,30,5) 
&= \begin{tabular}{|c|c|}
	\hline $E(28,30)^0_{4}$ & $P^1_{4}$ \\
	\hline $G(12,14, 4)^1$ & $G(12,14, 4)^0$ \\
	\hline
	\end{tabular} \\
G(12,14,4) 
&= \begin{tabular}{|c|c|}
	\hline $E(12,14)^0_{3}$ & $P^1_{3}$ \\
	\hline $G(4,6,   3)^1$ & $G(4,6,   3)^0$ \\
	\hline
	\end{tabular}
\end{align*}
The subsquares $P_3$ and $G(4,6,3)$ are base cases 
and can be looked up in the later appendices and the main part
of the paper.

\item Next, let $s = 6$, $\delta = 2^{6-2}=16$, and apply
\eqref{eqn:E} to get
\[
E(60,62)_{5} = 
	\begin{tabular}{|c|c|}
	\hline $E(28,30)^0_{4}$ & $L^1_{4}$ \\
	\hline $E(28,30)^1_{4}$ & $E(28,30)^0_{4}$ \\
	\hline
	\end{tabular}
\]
Next, let $s=5$, $\delta = 2^{s-2} = 2^3 = 8$. Then
\[
E(28,30)_{4} = 
		\begin{tabular}{|c|c|}
		\hline $E(12,14)^0_{3}$ & $L^1_{3}$ \\
		\hline $E(12,14)^1_{3}$ & $E(12,14)^0_{3}$ \\
		\hline
		\end{tabular}
\]
Lastly, let $s=4$, $\delta = 2^{s-2} = 2^2 = 4$. Then
\[
E(12,14)_{3} = 
		\begin{tabular}{|c|c|}
		\hline $E(4,6)^0_{2}$ & $L^1_{2}$ \\
		\hline $E(4,6)^1_{2}$ & $E(4,6)^0_{2}$ \\
		\hline
		\end{tabular}
\]
The subsquares $L_2$ and $E(4,6)_{2}$ are base cases 
and defined in the main section of the paper.
\end{itemize}
Here is 
$\mbox{gcs}(\alpha_{12,14} L_4)$:
\[
{\small
\begin{tabular}{|c|c|c|c|c|c|c|c||c|c|c|c|c|c|c|c|}
\hline 0 & 1  & -  & 3  & 4  & 5  & 6  & 7  & 8  & 9 & 10 & 11 & 12 & 13 & 14  & - \\
\hline 1 & 0 &  3  & 2  & 5  & 4  & 7  & 6  & 9  & 8 & 11 & 10 & 13 & 12  & -  & - \\
\hline 2 & 3 &  0  & 1  & 6  & 7  & 4  & 5 & 10 & 11  & 8  & 9 & 14  & - & 12  & - \\
\hline 3 & 2 &  1  & 0  & 7  & 6  & 5  & 4 & 11 & 10  & 9  & 8  & -  & -  & -  & - \\
\hline 4 & 5 &  -  & 7  & 0  & 1  & -  & 3 & 12 & 13 & 14  & -  & 8  & 9 & 10  & - \\
\hline 5 & 4 &  7  & 6  & 1  & 0  & 3  & 2 & 13 & 12  & -  & -  & 9  & 8  & -  & - \\
\hline 6 & 7 &  4  & 5  & 2  & 3  & 0  & 1 & 14  & - & 12  & - & 10  & -  & 8  & - \\
\hline 7 & 6 &  5  & 4  & 3  & 2  & 1  & 0  & -  & -  & -  & -  & -  & -  & -  & - \\
\hline \hline 8 & 9 &  - & 11 & 12 & 13 & 14  & -  & 0  & 1  & -  & 3  & 4  & 5  & 6  & - \\
\hline 9 & 8 & 11 & 10 & 13 & 12  & -  & -  & 1  & 0  & 3  & 2  & 5  & 4  & -  & - \\
\hline 10& 11&   8  & 9 & 14  & - & 12  & -  & 2  & 3  & 0  & 1  & 6  & -  & 4  & - \\
\hline 11& 10&   9  & 8  & -  & -  & -  & -  & 3  & 2  & 1  & 0  & -  & -  & -  & - \\
\hline 14&  -&  12  & - & 10  & -  & 8  & -  & 6  & -  & 4  & -  & 2  & -  & 0  & - \\
\hline - &12 & 15  & -  & -  & 8 & 11  & -  & -  & 4  & 7  & -  & -  & 0  & 3  & - \\
\hline 12& 13&   -  & -  & 8  & 9  & -  & -  & 4  & 5  & -  & -  & 0  & 1  & -  & - \\
\hline - & - &  -  & -  & -  & -  & -  & -  & -  & -  & -  & -  & -  & -  & -  & - \\
\hline
\end{tabular}
}
\]

\section{~}
\label{app:A}

\begin{align*}
\mbox{gcs}(\alpha_{4,6}L_3) 
& = \begin{tabular}{|c|c|c|c||c|c|c|c|}
  	\hline (0) & (1) &  2  & (3) & (4) & (5) & (6) &  7  \\
  	\hline (1) & (0) & (3) & (2) & (5) & (4) &  7  &  6  \\
  	\hline (2) & (3) & (0) & (1) & (6) &  7  & (4) &  5  \\
  	\hline (3) & (2) & (1) & (0) &  7  &  6  &  5  &  4  \\
	\hline
  	\hline (6) &  7  & (4) &  5  & (2) &  3  & (0) &  1  \\
  	\hline  5  & (4) & (7) &  6  &  1  & (0) & (3) &  2  \\
  	\hline (4) & (5) &  6  &  7  & (0) & (1) &  2  &  3  \\
  	\hline  7  &  6  &  5  &  4  &  3  &  2  &  1  &  0  \\
	\hline
	\end{tabular}
= 
	\begin{tabular}{|c||c|}
	\hline $E(4,6)_2$ & $\bullet$ \\
	\hline
	\hline $\bullet$ & $\bullet$ \\
	\hline
	\end{tabular}
\end{align*}
\begin{align*}
\mbox{gcs}(\alpha_{5,6}L_3) 
&= \begin{tabular}{|c|c|c|c||c|c|c|c|}
  	\hline (0) & (1) & (2) & (3) & (4) & (5) & (6) &  7  \\
  	\hline (1) & (0) & (3) & (2) & (5) & (4) &  7  &  6  \\
  	\hline (2) & (3) & (0) & (1) & (6) &  7  & (4) &  5  \\
  	\hline (3) & (2) & (1) & (0) &  7  &  6  &  5  &  4  \\
	\hline
  	\hline (4) & (5) & (6) &  7  & (0) & (1) & (2) &  3  \\
  	\hline (6) &  7  & (4) &  5  & (2) &  3  & (0) &  1  \\
  	\hline (5) & (4) &  7  &  6  & (1) & (0) &  3  &  2  \\
  	\hline  7  &  6  &  5  &  4  &  3  &  2  &  1  &  0  \\
	\hline
	\end{tabular}
= \begin{tabular}{|c||c|}
	\hline $E(5,6)_2$ & $\bullet$ \\
	\hline
	\hline $\bullet$ & $\bullet$ \\
	\hline
	\end{tabular}
\end{align*}
\begin{align*}
\mbox{gcs}(\alpha_{5,7}L_3) 
&= \begin{tabular}{|c|c|c|c||c|c|c|c|}
  	\hline (0) &  1  & (2) & (3) & (4) & (5) & (6) &  7  \\
  	\hline  1  &  0  & (3) & (2) & (5) & (4) &  7  &  6  \\
  	\hline (2) & (3) & (0) & (1) & (6) &  7  & (4) &  5  \\
  	\hline (3) & (2) & (1) & (0) &  7  &  6  &  5  &  4  \\
	\hline
  	\hline (4) &  5  & (6) &  7  & (0) &  1  & (2) &  3  \\
  	\hline  7  & (6) & (5) &  4  &  3  & (2) & (1) &  0  \\
  	\hline (6) & (7) &  4  &  5  & (2) & (3) &  0  &  1  \\
  	\hline  5  &  4  &  7  &  6  &  1  &  0  &  3  &  2  \\
	\hline
	\end{tabular}
= \begin{tabular}{|c||c|}
	\hline $E(5,7)_2$ & $\bullet$ \\
	\hline
	\hline $\bullet$ & $\bullet$ \\
	\hline
	\end{tabular}
\end{align*}
\begin{align*}
\mbox{gcs}(\alpha_{6,7}L_3) 
& = \begin{tabular}{|c|c|c|c||c|c|c|c|}
  	\hline  0  & (1) & (2) & (3) & (4) & (5) & (6) &  7  \\
  	\hline (1) & (0) & (3) & (2) & (5) & (4) &  7  &  6  \\
  	\hline  2  & (3) &  0  & (1) & (6) &  7  & (4) &  5  \\
  	\hline (3) & (2) & (1) & (0) &  7  &  6  &  5  &  4  \\
	\hline
  	\hline  4  & (5) & (6) &  7  &  0  & (1) & (2) &  3  \\
  	\hline (5) & (4) &  7  &  6  & (1) & (0) &  3  &  2  \\
  	\hline (7) &  6  & (5) &  4  & (3) &  2  & (1) &  0  \\
  	\hline  6  &  7  &  4  &  5  &  2  &  3  &  0  &  1  \\
	\hline
	\end{tabular}
= \begin{tabular}{|c||c|}
	\hline $E(6,7)_2$ & $\bullet$ \\
	\hline
	\hline $\bullet$ & $\bullet$ \\
	\hline
	\end{tabular}
\end{align*}

\section{~}
\label{app:B}
\begin{align*}
\mbox{gcs}(\alpha_{0,2}L_3) 
&=
	\begin{tabular}{|c|c|c|c||c|c|c|c|}
	\hline (2) & (3) & (0) & (1) & (6) &  7  & (4) &  5  \\
	\hline  1  & (0) & (3) & (2) &  5  & (4) & (7) &  6  \\
	\hline (0) & (1) & (2) & (3) & (4) & (5) &  6  &  7  \\
	\hline (3) & (2) & (1) & (0) &  7  &  6  &  5  &  4  \\
	\hline
	\hline (4) & (5) & (6) &  7  & (0) & (1) & (2) &  3  \\
	\hline (5) & (4) &  7  &  6  & (1) & (0) &  3  &  2  \\
	\hline (6) &  7  & (4) &  5  & (2) &  3  & (0) &  1  \\
	\hline  7  &  6  &  5  &  4  &  3  &  2  &  1  &  0  \\
	\hline
	\end{tabular}
= \begin{tabular}{|c||c|}
	\hline $A(0,2)_2$ & $\bullet$ \\
	\hline
	\hline $\bullet$ & $\bullet$ \\
	\hline
	\end{tabular}
\end{align*}
\begin{align*}
\mbox{gcs}(\alpha_{0,3}L_3) 
&= \begin{tabular}{|c|c|c|c||c|c|c|c|}
	\hline (3) &  2  &  1  & (0) & (7) & (6) & (5) &  4  \\
	\hline  1  & (0) &  3  & (2) & (5) &  4  & (7) &  6  \\
	\hline  2  &  3  & (0) & (1) & (6) & (7) &  4  &  5  \\
	\hline (0) & (1) & (2) & (3) &  4  &  5  &  6  &  7  \\
	\hline
	\hline (4) & (5) & (6) &  7  & (0) & (1) & (2) &  3  \\
	\hline (5) & (4) &  7  &  6  & (1) & (0) &  3  &  2  \\
	\hline (6) &  7  & (4) &  5  & (2) &  3  & (0) &  1  \\
	\hline  7  &  6  &  5  &  4  &  3  &  2  &  1  &  0  \\
	\hline
	\end{tabular}
= \begin{tabular}{|c||c|}
	\hline $A(0,3)_2$ & $\bullet$ \\
	\hline
	\hline $\bullet$ & $\bullet$ \\
	\hline
	\end{tabular}
\end{align*}
\begin{align*}
\mbox{gcs}(\alpha_{1,2}L_3) &= 
	\begin{tabular}{|c|c|c|c||c|c|c|c|}
	\hline (0) & (1) & (2) & (3) & (4) & (5) & (6) &  7  \\
	\hline (2) & (3) & (0) & (1) & (6) &  7  & (4) &  5  \\
	\hline (1) & (0) & (3) & (2) & (5) & (4) &  7  &  6  \\
	\hline (3) & (2) & (1) & (0) &  7  &  6  &  5  &  4  \\
	\hline
	\hline (4) & (5) & (6) &  7  & (0) & (1) & (2) &  3  \\
	\hline (5) & (4) &  7  &  6  & (1) & (0) &  3  &  2  \\
	\hline (6) &  7  & (4) &  5  & (2) &  3  & (0) &  1  \\
	\hline  7  &  6  &  5  &  4  &  3  &  2  &  1  &  0  \\
	\hline
	\end{tabular}
=
	\begin{tabular}{|c||c|}
	\hline $A(1,2)_2$ & $\bullet$ \\
	\hline
	\hline $\bullet$ & $\bullet$ \\
	\hline
	\end{tabular}
\end{align*}
\begin{align*}
\mbox{gcs}(\alpha_{1,3}L_3) &= 
	\begin{tabular}{|c|c|c|c||c|c|c|c|}
	\hline (0) & (1) & (2) & (3) & (4) &  5  & (6) &  7  \\
	\hline (3) &  2  &  1  & (0) &  7  & (6) & (5) &  4  \\
	\hline  2  &  3  & (0) & (1) & (6) & (7) &  4  &  5  \\
	\hline (1) & (0) & (3) & (2) &  5  &  4  &  7  &  6  \\
	\hline
	\hline (4) & (5) & (6) &  7  & (0) & (1) & (2) &  3  \\
	\hline (5) & (4) &  7  &  6  & (1) & (0) &  3  &  2  \\
	\hline (6) &  7  & (4) &  5  & (2) &  3  & (0) &  1  \\
	\hline  7  &  6  &  5  &  4  &  3  &  2  &  1  &  0  \\
	\hline
	\end{tabular}
=
	\begin{tabular}{|c||c|}
	\hline $A(1,3)_2$ & $\bullet$ \\
	\hline
	\hline $\bullet$ & $\bullet$ \\
	\hline
	\end{tabular}
\end{align*}
\begin{align*}
\mbox{gcs}(\alpha_{2,3}L_3) &= 
	\begin{tabular}{|c|c|c|c||c|c|c|c|}
	\hline (0) & (1) & (2) & (3) &  4  & (5) & (6) &  7  \\
	\hline (1) & (0) & (3) & (2) & (5) & (4) &  7  &  6  \\
	\hline  3  & (2) &  1  & (0) & (7) &  6  & (5) &  4  \\
	\hline (2) & (3) & (0) & (1) &  6  &  7  &  4  &  5  \\
	\hline
	\hline (4) & (5) & (6) &  7  & (0) & (1) & (2) &  3  \\
	\hline (5) & (4) &  7  &  6  & (1) & (0) &  3  &  2  \\
	\hline (6) &  7  & (4) &  5  & (2) &  3  & (0) &  1  \\
	\hline  7  &  6  &  5  &  4  &  3  &  2  &  1  &  0  \\
	\hline
	\end{tabular}
=
	\begin{tabular}{|c||c|}
	\hline $A(2,3)_2$ & $\bullet$ \\
	\hline
	\hline $\bullet$ & $\bullet$ \\
	\hline
	\end{tabular}
\end{align*}

\section{~}
\label{app:0}

\begin{align*}
U(5,6) &= 
	\begin{tabular}{|c|c|c|c||c|c|c|c|}
	\hline (0) & (1) & (2) & (3) & (4) & (5) & (6) &  7  \\
	\hline (1) & (0) & (3) & (2) & (5) & (4) &  7  &  6  \\
	\hline (2) & (3) & (0) & (1) & (6) &  7  & (4) &  5  \\
	\hline (3) & (2) & (1) & (0) &  7  &  6  &  5  &  4  \\
	\hline
	\hline (4) & (5) & (6) & (7) & (0) & (1) & (2) &  3  \\
	\hline (5) & (4) & (7) & (6) & (1) & (0) &  3  &  2  \\
	\hline (6) & (7) & (4) & (5) & (2) &  3  & (0) &  1  \\
	\hline (7) & (6) & (5) & (4) &  3  &  2  &  1  &  0  \\
	\hline 
	\end{tabular} \\
~ \\
U(5,7) &=
	\begin{tabular}{|c|c|c|c||c|c|c|c|}
	\hline (0) &  1  & (2) & (3) & (4) & (5) & (6) &  7  \\
	\hline  1  &  0  & (3) & (2) & (5) & (4) &  7  &  6  \\
	\hline (2) & (3) & (0) & (1) & (6) &  7  & (4) &  5  \\
	\hline (3) & (2) & (1) & (0) &  7  &  6  &  5  &  4  \\
	\hline
	\hline (4) &  5  & (6) & (7) & (0) & (1) & (2) &  3  \\
	\hline  5  &  4  & (7) & (6) & (1) & (0) &  3  &  2  \\
	\hline (6) & (7) & (4) & (5) & (2) &  3  & (0) &  1  \\
	\hline (7) & (6) & (5) & (4) &  3  &  2  &  1  &  0  \\
	\hline 
	\end{tabular} \\
~ \\
U(6,7) &=
	\begin{tabular}{|c|c|c|c||c|c|c|c|}
	\hline  0  & (1) & (2) & (3) & (4) & (5) & (6) &  7  \\
	\hline (1) & (0) & (3) & (2) & (5) & (4) &  7  &  6  \\
	\hline  2  & (3) &  0  & (1) & (6) &  7  & (4) &  5  \\
	\hline (3) & (2) & (1) & (0) &  7  &  6  &  5  &  4  \\
	\hline
	\hline  4  & (5) & (6) & (7) & (0) & (1) & (2) &  3  \\
	\hline (5) & (4) & (7) & (6) & (1) & (0) &  3  &  2  \\
	\hline  6  & (7) &  4  & (5) & (2) &  3  & (0) &  1  \\
	\hline (7) & (6) & (5) & (4) &  3  &  2  &  1  &  0  \\
	\hline 
	\end{tabular}
\end{align*}

\section{~}
\label{app:1}

\begin{align*}
V(4,5) = V(4,6) &=
	\begin{tabular}{|c|c|c|c||c|c|c|c|}
	\hline (0) & (1) &  2  & (3) & (4) & (5) & (6) &  7  \\
	\hline (1) & (0) & (3) & (2) & (5) & (4) &  7  &  6  \\
	\hline (2) & (3) & (0) & (1) & (6) &  7  & (4) &  5  \\
	\hline (3) & (2) & (1) & (0) &  7  &  6  &  5  &  4  \\
	\hline
	\hline (4) & (5) & (6) &  7  & (0) & (1) &  2  & (3) \\
	\hline (5) & (4) &  7  &  6  & (1) & (0) & (3) & (2) \\
	\hline (6) &  7  & (4) &  5  & (2) & (3) & (0) & (1) \\
	\hline  7  &  6  &  5  &  4  & (3) & (2) & (1) & (0) \\
	\hline 
	\end{tabular} \\
~ \\
V(5,6) &= 
	\begin{tabular}{|c|c|c|c||c|c|c|c|}
	\hline (0) & (1) & (2) & (3) & (4) & (5) & (6) &  7  \\
	\hline (1) & (0) & (3) & (2) & (5) & (4) &  7  &  6  \\
	\hline (2) & (3) & (0) & (1) & (6) &  7  & (4) &  5  \\
	\hline (3) & (2) & (1) & (0) &  7  &  6  &  5  &  4  \\
	\hline
	\hline (4) & (5) & (6) &  7  & (0) & (1) & (2) & (3) \\
	\hline (5) & (4) &  7  &  6  & (1) & (0) & (3) & (2) \\
	\hline (6) &  7  & (4) &  5  & (2) & (3) & (0) & (1) \\
	\hline  7  &  6  &  5  &  4  & (3) & (2) & (1) & (0) \\
	\hline 
	\end{tabular} \\
~ \\
V(5,7) &= 
	\begin{tabular}{|c|c|c|c||c|c|c|c|}
	\hline (0) &  1  & (2) & (3) & (4) & (5) & (6) &  7  \\
	\hline  1  &  0  & (3) & (2) & (5) & (4) &  7  &  6  \\
	\hline (2) & (3) & (0) & (1) & (6) &  7  & (4) &  5  \\
	\hline (3) & (2) & (1) & (0) &  7  &  6  &  5  &  4  \\
	\hline
	\hline (4) & (5) & (6) &  7  & (0) &  1  & (2) & (3) \\
	\hline (5) & (4) &  7  &  6  &  1  &  0  & (3) & (2) \\
	\hline (6) &  7  & (4) &  5  & (2) & (3) & (0) & (1) \\
	\hline  7  &  6  &  5  &  4  & (3) & (2) & (1) & (0) \\
	\hline 
	\end{tabular} \\
~ \\
V(6,7) &=
	\begin{tabular}{|c|c|c|c||c|c|c|c|}
	\hline  0  & (1) & (2) & (3) & (4) & (5) & (6) &  7  \\
	\hline (1) & (0) & (3) & (2) & (5) & (4) &  7  &  6  \\
	\hline  2  & (3) &  0  & (1) & (6) &  7  & (4) &  5  \\
	\hline (3) & (2) & (1) & (0) &  7  &  6  &  5  &  4  \\
	\hline
	\hline (4) & (5) & (6) &  7  &  0  & (1) & (2) & (3) \\
	\hline (5) & (4) &  7  &  6  & (1) & (0) & (3) & (2) \\
	\hline (6) &  7  & (4) &  5  &  2  & (3) &  0  & (1) \\
	\hline  7  &  6  &  5  &  4  & (3) & (2) & (1) & (0) \\
	\hline 
	\end{tabular}
\end{align*}
\end{document}